\documentclass{article}

\usepackage{amsfonts}
\usepackage{pb-diagram}
\usepackage{amsmath,amssymb,amsthm}
\usepackage{multicol}
\usepackage{color}
\usepackage{hyperref}
\usepackage{graphicx}
\usepackage[utf8x]{inputenc}
\usepackage[english]{babel}

\newtheorem{Def}{Definition}[section]
\newtheorem{theorem}{Theorem}

\newtheorem{Prop}[Def]{Property}

\newcommand{\pia}{\pi_{\mathfrak{a}}}

\newcommand{\gf}{\mathfrak{g}}

\newcommand{\af}{\mathfrak{a}}

\newcommand{\bff}{\mathfrak{b}}
\newcommand{\afb}{\mathfrak{a}_{\bot}}
\newcommand{\hf}{\mathfrak{h}}
\newcommand{\hfg}{\hf_{\gf}}
\newcommand{\hfa}{\hf_{\af}}

\newcommand{\sfr}{\mathfrak{s}}

\begin{document}

\title{Splints of root systems for special Lie subalgebras} 

\author{V.D.~Lyakhovsky$^1$, A.A.~Nazarov$^{1,2}$, P.I.~Kakin$^{1}$\\
  {\small $^1$ Department of High-energy and elementary particle physics,}\\ {\small St Petersburg State University}\\
  {\small 198904, Saint-Petersburg, Russia,}\\
  {\small $^{2}$ e-mail: antonnaz@gmail.com}}
\date{}
\maketitle

\begin{abstract}
  Splint is a decomposition of root system into union of root systems. Splint of root system for
  simple Lie algebra appears naturally in studies of (regular) embeddings of reductive subalgebras.
  Splint can be used to construct branching rules. We consider special embedding of Lie subalgebra
  to Lie algebra. We classify projections of algebra root systems using extended Dynkin diagrams and
  single out the conditions of splint appearance and coincidence of branching coefficients with
  weight multiplicities. While such a coincidence is not very common it is connected with Gelfand-Tsetlin basis. 
\end{abstract}

\section{Introduction}
\label{sec:introduction}

The notion of splint was introduced by David Richter in the paper \cite{richter2008splints}. Splint
is the decomposition of root system into disjoint union of images of two or more embeddings of some
other root systems. Embedding $\phi$ of a root system $\Delta_1$ into a root system $\Delta$ is a
bijective map of roots of $\Delta_{1}$ to a (proper) subset of $\Delta$ that commutes with vector
composition law in $\Delta_{1}$ and $\Delta$.
\begin{equation*}
\phi:\Delta_1 \longrightarrow \Delta
\end{equation*}
\begin{equation*}
\phi \circ (\alpha + \beta) =\phi \circ \alpha + \phi \circ \beta,
\,\,\, \alpha,\beta \in \Delta_1
\end{equation*}

Note that the image $Im(\phi)$ is not required to inherit the root system properties except the
addition rules equivalent to the addition rules in $\Delta_{1}$ (for pre-images). If an embedding
$\phi$ preserves the angles between the roots it is called ``metric''. Two embeddings $\phi_1$ and
$\phi_2$ can splinter $\Delta$ when the latter can be presented as a disjoint union of images
$Im(\phi_1)$ and $Im(\phi_2)$.

Why one would consider non angle-preserving maps of root systems? Additive properties of root system
determine the structure of Verma module:
    \begin{equation*}
      M^{\mu}=U(\gf)\underset{U(\bff_{+})}{\otimes} D^{\mu}(\bff_{+})=\{(E^{-\alpha_{1}})^{n_{1}}\dots (E^{-\alpha_{s}})^{n_{s}} \left|v_{\mu}\right>\}_{\alpha_{i}\in\Delta^{+}}^{n_{i}=0,1,\dots}
    \end{equation*}
Here $\bff_{+}$ is the Borel subalgebra of $\gf$, $D^{\mu}$ its one-dimensional representation,
$\left|v_{\mu}\right>$ is the highest weight vector and $E^{-\alpha_{j}}$ are lowering operators
that are in correspondence with positive roots $\alpha_{j}\in \Delta^{+}$. 
Weyl character formula expresses character of irreducible representation as a combination of
characters of Verma modules: 
\begin{equation*}
  \mathrm{ch} L^{\mu}=\frac{\sum_{w\in W} \varepsilon(w) e^{w(\mu+\rho)-\rho}}{\sum_{w\in W}\varepsilon(w) e^{w\rho-\rho}}=\sum_{w\in W} \varepsilon(w)\; \mathrm{ch} M^{w(\mu+\rho)-\rho}
\end{equation*}
Here the sum is over elements of the Weyl group and their actions $w\triangleright \mu
=w(\mu+\rho)-\rho$ depend on angles between roots. But we can rewrite Weyl character formula
representing Weyl group elements as the products of reflections $s_{\alpha}, \alpha \in S$ in
hyperplanes orthogonal to simple roots: $w=s_{\alpha_{1}}\cdot s_{\alpha_{2}}\dots$. Reflections act
on root system by permutations, so the composition of such an action with embedding $\phi$ is easily
obtained. The orbit of Weyl groups action $w\triangleright \mu, w\in W$, can be constructed by
subtractions of roots from the highest weight $\mu$. Let's denote Dynkin labels of $\mu$ by
$(\mu_{1},\dots \mu_{r})$. Then $\mu-\mu_{1}\alpha_{1}, \mu-\mu_{2}\alpha_{2},\dots,
\mu-\mu_{r}\alpha_{r}, \alpha_{i}\in S$, are the points of the orbit that can be obtained by
elementary reflections $s_{\alpha_{i}}, \alpha_{i}\in S$. Next set of points that are obtained by
two consecutive reflections $s_{\alpha_{i}}\cdot s_{\alpha_{j}}\triangleright\mu$ are obtained as
the subtraction $\mu-\mu_{i}\alpha_{i}-\mu_{j} (s_{\alpha_{i}}\alpha_{j})$. Continuing this way we
construct the image of singular element $\Psi=\sum_{w\in W} \varepsilon(w) e^{w(\mu+\rho)-\rho}$
after the embedding $\phi$ \cite{2011arXiv1111.6787L}. The multiplicities in the character of Verma
module are unchanged by the embedding since they are determined by additive properties of the roots.
So we see that weight multiplicities of irreducible modules are also preserved by the embedding.

Root system of regular subalgebra is contained in the root system of algebra so splint is useful in
computation of branching coefficients. 

In the paper \cite{2011arXiv1111.6787L} it was proven that the existence of splint leads to the
coincidence of branching coefficients with weight multiplicities under certain conditions. We denote
Lie algebra by $\gf$ and consider it's subalgebra $\af$. If $\af$ is a regular subalgebra, its root
system $\Delta_{\af}$ is contained in $\Delta_{\gf}$. Branching coefficients $b^{(\mu)}_{\nu}$
appear in decomposition of irreducible representation $L^{\mu}_{\gf}$ of $\gf$ to the sum of
irreducible representations of $\af$:
\begin{equation}
  \label{eq:1}
  L^{\mu}_{\gf}=\bigoplus_{\nu} b^{(\mu)}_{\nu} L^{\nu}_{\af}
\end{equation}

Assume that root system $\Delta_{\gf}$ splinters as $\Delta_{\gf}=\Delta_{\af} \cup
\phi(\Delta_{\sfr})$, where $\phi$ is an embedding of root system $\Delta_{\sfr}$ of some semisimple
Lie algebra $\sfr$. Then branching coefficients $b^{(\mu)}_{\nu}$ for the reduction
$L^{(\mu)}_{\gf\downarrow \af}$ coincide with weight multiplicities $m^{(\tilde \mu)}_{\nu}$ in
$\sfr$-representations provided certain technical condition holds \cite{2011arXiv1111.6787L}.
Highest weight $\tilde\mu$ of $\sfr$-representation is calculated from Dynkin labels of highest
weight $\mu$ of $\gf$-module.

In the next section \ref{sec:splints-root-systems} we review the classification of splints for simple
root systems and results for branching coefficients. Then we move to the study of special embeddings
which is the main subject of the present paper. We consider special embeddings of Lie subalgebras
into a Lie algebra. In this case root system of the subalgebra is not contained in the root system of
the algebra. So the original motivation for splints is not applicable. But we can consider the
projection of root system of the algebra on the root space of the subalgebra. Such a projection is not a
root system anymore, but it satisfies milder conditions (Section \ref{sec:spec-embedd-proj}). It is
possible to classify most projections using Dynkin diagrams augmented with multiplicities. We then
define splint for such 'weak' root systems and state the conditions of its appearance and
implications for the calculation of branching coefficients (Section \ref{sec:splints-spec-embedd}).

Use of representation theory of the subalgebra allows us to classify all the splints for the projections
of algebra root system. We also apply this method to metric splints and regular subalgebras and get
a unified treatment. We obtain Gelfand-Tsetlin rules for regular and special embeddings this way.

In conclusion \ref{sec:conclusion} we discuss the cases when the projection of the root system does
not fall into the classification mentioned above.

\section{Splints of root systems and regular subalgebras}
\label{sec:splints-root-systems}
The notion and classification of
splints for simple root systems were introduced in the paper \cite{richter2008splints}.

Consider a simple Lie algebra $\mathfrak{g}$ and its regular subalgebra $%
\mathfrak{a}\hookrightarrow \mathfrak{g}$ such that $\mathfrak{a}$
is a
reductive subalgebra $\mathfrak{a \subset g}$ with correlated root spaces: $%
\mathfrak{h}_{\mathfrak{a}}^{\ast }\subset \mathfrak{h}_{\mathfrak{g }%
}^{\ast }$. Let $\mathfrak{a}^{\mathfrak{s}}$ be a semisimple summand of
$\mathfrak{a}$,
this means that $\mathfrak{a}=\mathfrak{a}^{\mathfrak{s}} \oplus \mathfrak{u}(1)\oplus %
\mathfrak{u}(1)\oplus \dots$. We shall consider $\mathfrak{a}^{\mathfrak{s}}$
to be a proper regular subalgebra and $\mathfrak{a}$ to be the
maximal subalgebra with $\mathfrak{a}^{\mathfrak{s}}$ fixed that is the rank
$r$ of $\frak{a}$ is equal to that of $\mathfrak{g}$.

Computation of branching coefficients relies on roots $\Delta_{\gf}\setminus \Delta_{\af}$
\cite{2010arXiv1007.0318L} so splint  is naturally connected to branching coefficients for regular
subalgebra.

Only three types of splints are
injective and thus are naturally connected to branching \cite{richter2008splints}:
\begin{equation}
\label{eq:1}
\begin{array}{cc||c|c}
\hbox{type} & \hspace{0.25in}\Delta \hspace{0.25in} & \hspace{0.25in}\Delta
_{\frak{a}}\hspace{0.25in} & \hspace{0.25in}\Delta _{\sfr}\hspace{0.25in}
\\ \hline\hline
\hbox{(i)} & G_{2} & A_{2} & A_{2} \\
& F_{4} & D_{4} & D_{4} \\ \hline
\hbox{(ii)} & B_{r}(r\geq 2) & D_{r} & \oplus ^{r}A_{1} \\
(*)& C_{r}(r\geq 3) & \oplus ^{r}A_{1} &  D_{r} \\ \hline
\hbox{(iii)} & A_{r}(r\geq 2) & A_{r-1}\oplus u\left( 1\right)  & \oplus
^{r}A_{1} \\
& B_{2} & A_{1}\oplus u\left( 1\right)  & A_{2}
\end{array}
\end{equation}

Each row in the table gives a splint $(\Delta _{\frak{a}},\Delta _{\sfr})
$ of the simple root system $\Delta $. In the first two types both $\Delta _{%
\frak{a}}$ and $\Delta _{\sfr}$ are embedded metrically. Stems in the
first type splints are equivalent and in the second are not. In the third
type splints only $\Delta _{\frak{a}}$ is embedded metrically. The summands $%
u\left( 1\right) $ are added to keep $r_{\frak{a}}=r$. This does not change
the principle properties of branching but makes it possible to use the
multiplicities of $\sfr$ -modules without further projecting their
weights.

Note that in the case of $C_{r}$-series (marked with the star in table \eqref{eq:1}) metrical
embedding $\Delta_{D_{r}}\to \Delta_{C_{r}}$ does not lead to the appearance of regular subalgebra
$\mathfrak{a}$ (See \cite{dynkin1952semisimple}).

In the paper \cite{2011arXiv1102.1702L} it was shown that

\begin{Prop}
There is a decomposition of a singular element of the algebra $\gf$ into singular elements of the subalgebra $\sfr$:
\begin{equation}
 \Phi ^{\mu }_{\gf} = \sum_{w\in W_{\frak{a}}}\epsilon\left( w\right)w\left( e^{\mu +\rho _{\gf}}\phi\left( e^{-\widetilde{\mu }}\Psi^{\widetilde{\mu }}_{\sfr}\right)\right),
\end{equation}
where $\epsilon\left( w\right)$ is a determinant of the element $w$ of the Weyl group $W_{\af}$ of
the subalgebra $\af$. The action of Weyl group on the weights is extended to the algebra of formal
exponents by the rule $w\left(e^{\nu}\right)=e^{w\nu}, w\in W$. Hence, the multiplicity $M_{\left(
      \sfr\right) \widetilde{\nu }}^{\widetilde{\mu }}$ of the weight $\widetilde{\nu }$ from the
  weight diagram of the algebra $\sfr$ with the highest weight $\tilde\mu$ defines the branching
  coefficient $b_{\nu }^{(\mu )}$ for the highest weight $\nu =\left( \mu -\phi \left(
      \widetilde{\mu }-\widetilde{\nu }\right) \right) $:
\begin{equation}
b_{\left( \mu -\phi \left( \widetilde{\mu }-\widetilde{\nu }\right) \right)
}^{(\mu )}=M_{\left( \sfr\right) \widetilde{\nu }}^{\widetilde{\mu }}. 
\label{bran1}
\end{equation}
\end{Prop} 

\section{Special embeddings and projections of root system}
\label{sec:spec-embedd-proj}

The study of specialy embedded subalgebras traces back to fundamental papers by Eugene Dynkin
\cite{dynkin1952semisimple,dynkin1952maximal} where he called such subalgebras ``S-subalgebras'' to
distinguish from regular or ``R-subalgebras'' that have root systems obtained by dropping some roots
of the algebra root system. Thus regular subalgebras are easy to construct by dropping nodes from
extended Dynkin diagram of the algebra, while the case of special subalgebras is more difficult. Complete
classification for exceptional Lie algebras was obtained recently in \cite{minchenko2006semisimple}.
The algorithm of the special subalgebras construction is available in GAP package but still requires
manual intervention \cite{de2011constructing}.  

Assume that Lie algebra $\gf$ is simple. Let's denote a subalgebra by
$\af$. We denote corresponding Cartan subalgebras by $\hfg$ and $\hfa$ and identify them with the dual
spaces $\hfg^{*}$, $\hfa^{*}$ using Killing forms of $\gf$ and $\af$.

To construct an embedding $\af\to\gf$ consider some representation $L^{\nu}_{\af}$ as a subspace of
Lie algebra $\gf$. Then one needs to check that the generators of $\af$ can be presented as
linear combinations of generators of $\gf$. We can identify Cartan subalgebra $\hfa\subset \af$ with
dual space $\hfa^{*}$ using Killing form. Then the root system $\Delta_{\gf}$ of $\gf$ can be
projected to $\hfa^{*}$ using the expression of $\hfa$-generators through $\hfg$-generators.


This projection is described in the classical papers \cite{dynkin1952semisimple,dynkin1952maximal} by the following theorem:

\begin{theorem}\label{dyn0}
  If a representation $L^{\mu}_{\gf}$ of the algebra $\gf$ induces a representation
  $L^{\tilde\mu}_{\af}$ on the subalgebra $\af$ then the weight system of $L^{\tilde\mu}_{\af}$ can
  be obtained from $L^{\mu}_{\gf}$ by an orthogonal projection of $\hfg^{*}$ on $\hfa^{*}$.

  \cite{dynkin1952maximal}. 
\end{theorem} 

In relation to the adjoint representation of $\gf$ this theorem means that the projection of the root system $\Delta_{\gf}$ to $\hfa^{*}$ is a weight system of some
finite-dimensional but not irreducible representation of $\af$. Moreover this reducible representation contains the adjoint representation of $\af$. And if we denote such a projection by
\begin{equation}
  \label{eq:2}
  \Delta'=\pia\left(\Delta_{\gf}\right),
\end{equation}
then the roots of $\af$ will be contained in $\Delta'$. Thus the system $\Delta'$ can be a root
system but it might not be reduced and might contain some vectors with multiplicities greater than
one.

The following theorems \cite{dynkin1952semisimple} also concern the properties of $\Delta'$. 

\begin{theorem}\label{dyn1}
  Every special subalgebra $\af$ of a semisimple Lie algebra $\gf$ is integer, i.e. the projections of
  roots of $\gf$ to $\hfa^{*}$ are linear combinations of simple roots of $\af$ with integer
  coefficients \cite{dynkin1952semisimple}
\end{theorem}

\begin{theorem}\label{dyn2}
  If $\af$  is a semisimple subalgebra of a semisimple Lie algebra $\gf$ and the generator $e_{\alpha}$ corresponding to
  the root $\alpha$ of $\af$ is presented as a linear combination $e_{\alpha}=\sum_{\beta}
  e_{\beta}$ of $\gf$-generators
  corresponding to the roots $\beta$ of $\gf$, then the roots $\beta$ are projected into $\alpha$,
  $\pia(\beta)=\alpha$. 
  \cite{dynkin1972semisimple,dynkin1952semisimple}. 
\end{theorem}

From these theorems we see that the root system $\Delta_{\af}$ is metrically embedded (in the sense
of splint) into the projection $\Delta'$. Moreover, the multiplicities of some roots
$\alpha\in\Delta_{\af}$ in this projection $\Delta'$ are greater than one. 

\begin{figure}[h!bt]
  \noindent\centering{
   \includegraphics[width=80mm]{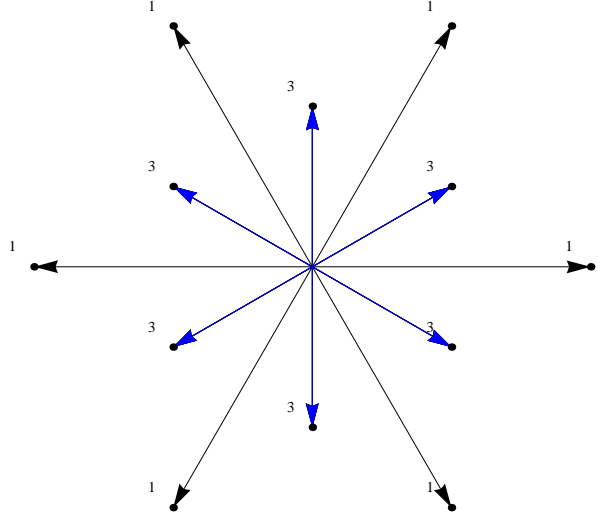}
  }
  \caption{Projection of $D_{4}$ ($so(8)$)-root system onto the root space of the special subalgebra $A_{2}$
    ($su(3)$). Note that this roots system is $G_{2}$ but with non-trivial multiplicities. }
 \label{fig:d4-a2_splint}
\end{figure}

As an example we present a projection $\Delta'$ of the root system of $D_{4}$ to the root system of the special
subalgebra $A_{2}$ (Fig. \ref{fig:d4-a2_splint}).

\section{Splints for special embeddings}
\label{sec:splints-spec-embedd}


We want to classify all the cases when branching coefficients coincide with weight multiplicities of
some other algebra representations. Such a coincidence is possible when the projection $\Delta'$ of
$\gf$ root system admits a splint $\Delta'=\varphi_{\af}(\Delta_{\af})\cup
\varphi_{\sfr}(\Delta_{\sfr})$, where $\varphi_{\af}$ and $\varphi_{\sfr}$ are the embeddings of
corresponding root systems. The embedding $\varphi_{\af}$ is metric and trivial. Moreover, the
projection of singular element $\pia\left(\Psi^{\mu}_{\gf}\right)$ should admit a decomposition into
linear combination of singular elements of irreducible representations of $\sfr$:
\begin{equation}
  \label{eq:4}
  \pia\left(\Psi^{\mu}_{\gf}\right)=\sum_{\nu} \varkappa_{\nu}e^{\nu}\Psi^{\tilde\mu}_{\sfr},
\end{equation}
with integer coefficients $\varkappa_{\nu}$, where the sum is over the set of weights $\nu$ that
will be determined later. First we classify all the splints of the projection of a root system into
a union of simple root systems for special embeddings, and then discuss the decomposition of singular
elements.

Projection of the $\gf$ root system $\Delta'$ coincides with the projection of the weight diagram of the
adjoint representation of algebra $\gf$ with the exception of zero weight. The adjoint representation of
$\gf$ contains adjoint representation of $\af$ and can be decomposed as
\begin{equation}
  \label{eq:3}
  \mathrm{ad}_{\gf}=\mathrm{ad}_{\af}\oplus M_{\af}^{\chi},
\end{equation}
where $M^{\chi}_{\af}$ is not necessary irreducible representation of $\af$ called
``characteristic'' in \cite{dynkin1952semisimple}. 

The weight diagram of $M^{\chi}_{\af}$ coincides with the root system $\Delta_{\sfr}$ with the exception of
zero weight. So we need to find all the representations for all simple Lie algebras $\af$, such that
their weight diagrams are root systems after the exclusion of zero weight. In order to do so we must
note, that all the weights of $M^{\chi}_{\af}$ should have length not less than length of shortest
root of $\af$, since $\af$ is an integer subalgebra of $\gf$ \ref{dyn1}. Moreover, $M^{\chi}_{\af}$
should contain weights of at most two different lengths. 

If the root system $\Delta_{\gf}$ contains roots $\alpha: \alpha\perp \beta,\; \forall \beta\in
\Delta_{\af},$ orthogonal to the root system of the subalgebra, one cannot proceed with the decomposition
of singular element \eqref{eq:4}. The elements $e^{\nu}$ of $\pia\left(\Psi^{\mu}_{\gf}\right)$
should be augmented with the dimensions of representations of algebra $\afb$ spanned by generators
corresponding to orthogonal roots \cite{2010arXiv1007.0318L}. Then there will be non-trivial
multiplicities in the right hand side of \eqref{eq:4} and no coincidence of branching coefficients
with weight multiplicities of $\sfr$-representations. One can write more cumbersome relation between
multiplicities and branching coefficients, but it is out of scope for the present paper.

The multiplicity of zero weight in the adjoint representation is equal to the rank of the algebra. So if
$\mathrm{rank}\gf-\mathrm{rank}\af=1$, the representation in question $M^{\chi}_{\af}$ must be
multiplicity-free. If $\mathrm{rank}\gf-\mathrm{rank}\af>1$ only zero weight can have non-trivial
multiplicities. 

The simplest class of multiplicity-free representations is called ``strongly multiplicity free''
\cite{lehrer2006strongly} and consists of multiplicity free representations with weight systems
admitting strict ordering $\nu_{1}<\nu_{2} \Leftrightarrow \nu_{1}=\nu_{2}+n \alpha$, where $n\in
\mathbb{N}, \alpha\in \Delta^{+}_{\af}$ and for all $\nu_{1},\nu_{2}$ either $\nu_{1}<\nu_{2}$ or
$\nu_{2}<\nu_{1}$.

The list of strongly multiplicity free representations consists of (first) fundamental representations
for series $A_{r}, B_{r}, C_{r}$, exceptional Lie algebra $G_{2}$ (7-dimensional representation) and
all $A_{1}$ representations. 

Fundamental weights of the algebra $A_{r}$ are shorter than its roots, but the projections of $\gf$
weights must be given by linear combinations of the roots of the  special subalgebra $\af$ with integer
coefficients, so the first class of strongly multiplicity free representations does not produce such a
projection. 

Nevertheless, the union of diagrams of two fundamental representations of $A_{2}$ with the root
system of $A_{2}$ produces weight diagram of $G_{2}$. This case corresponds to splint
$\Delta_{G_{2}}=\varphi_{1}( \Delta_{A_{2}})\cup \varphi_{2}(\Delta_{A_{2}})$ which is connected
with regular subalgebra $A_{2}\subset G_{2}$ (see Section \ref{sec:splints-root-systems}).

First fundamental representation of $B_{r}$ immediately gives us the splint $\Delta'=\pi_{B_{r}}\left(
\Delta_{D_{r+1}}\right) = \Delta_{B_{r}}\cup \Delta_{A_{1}+\dots+A_{1}}$ corresponding to Gelfand-Tsetlin
multiplicity-free branching for the special embedding $\mathfrak{so}(2r+1)\to \mathfrak{so}(2r+2)$. 

The length of the first fundamental weight of  $C_{r}$ is less than the length of its short root, so
this case can not correspond to an integer subalgebra. 

If there are vectors of different length in the projection $\Delta'$ on the subalgebra $A_{1}$, such a
projection is not multiplicity free, and $\Delta_{\sfr}$ contains parallel roots. The simplest case
is the special embedding $A_{1}\to A_{2}$ with index 4 where $\Delta_{\sfr}$ is the root system
$BC_{1}$. Such systems do not correspond to semisimple Lie algebras, so we can not speak of a
coincidence of branching coefficients for the reduction $\gf\downarrow \af$ with weight
multiplicities of $\sfr$-representations.

The weight diagram of the seven-dimensional representation of $G_{2}$ together with $G_{2}$ root system
form the projection of the root system $B_{3}$ which corresponds to special embedding 
$G_{2}\to B_{3}$. Branching coefficients in this case coincide with weight multiplicities of
representations of algebra $\sfr=A_{2}$ with image of root system consisting of short roots of
$G_{2}$. 

The complete classification of multiplicity-free irreducible representations was obtained in
\cite{howe1995perspectives} (see also \cite{stembridge2003multiplicity}). It consists of minuscule
and quasi-minuscule representations.

The weight $\mu$ is minuscule if $\left<\mu,\alpha^{\vee}\right>\leq 1$ for all $\alpha\in
\Delta^{+}$ and all the weights of the irreducible representation $L^{\mu}$ lie on the Weyl
group orbit of $\mu$. A weight is quasi-minuscule if $\left<\mu,\alpha^{\vee}\right>\leq 2$ for all
$\alpha\in \Delta^{+}$ and all non-zero weights lie on the single Weyl group orbit.

Minuscule representations are well-studied and very useful in computations of e.g. tensor products
\cite{stembridge2003multiplicity,stembridge2001computational}. The minuscule representations are
indexed by the weight lattice modulo the root lattice. There is a unique quasi-minuscule
representation that is not minuscule for each simple Lie algebra. The multiplicity of the zero 
weight in quasi-minuscule representation is the number of short nodes of the Dynkin diagram.

List of minuscule representations consists of tensor powers of vector
representation for the series $A_{r}$; spin representations for the series
$B_{r}$; vector representations for $C_{r}$; vector and two half-spin for $D_{r}$; two
27-dimensional representations for $E_{6}$ and 56-dimensional representation of $E_{7}$. 

Quasi-minuscule representations that are not minuscule are: adjoint representation of $A_{r}$,
vector representation for $B_{r}$, $2r^{2}-r-1$-dimensional representation for $C_{r}$, adjoint
representations for $D_{r}, E_{6}, E_{7}, E_{8}$, 26-dimensional representation for $F_{4}$ and
7-dimensional of $G_{2}$.


Classification of multiplicity-free irreducible highest weight modules $L^{\mu}$ in
\cite{howe1995perspectives,stembridge2003multiplicity} consists of the following classes:
\begin{itemize}
\item (1) $\mu$ is minuscule,
\item (2) $\mu$ is quasi-minuscule and Lie algebra has only one short simple root,
\item (3) Lie algebra is $\mathfrak{sp}(6)$ and $\mu=\omega_{1}$, or
\item (4) Lie algebra is $\mathfrak{sl}(n + 1)$ and  $\mu= m\omega_{1}$ or $\mu  = m\omega_{n}$ 
\end{itemize}

We see that case (1) contains strongly multiplicity free modules of series $A_{r}$ and exceptional
Lie algebra $G_{2}$, class (2) contains strongly multiplicity free modules of $ B_{r}$ and $C_{r}$ ,
strongly multiplicity free modules of $A_{1}$ are included in class (4).

We do not need to consider minuscule representations for simply-laced Lie algebras, since in this
case length of the weights of representation $M^{\chi}_{\af}$ is less than length of $\af$ roots,
that contradicts to $\af$ being an integer subalgebra. For the series $B_{r}$ we have already
discussed minuscule spin representation, since it is strongly multiplicity free. The last case in
(1) is the vector representation for $C_{r}$. It is $2r$-dimensional, but the length of its highest weight $\omega_{1}$ is less
than the length of the short root of algebra $C_{r}$. So it can not be a characteristic representation
for the special subalgebra. 

Case (2) consists of the series $B_{r}$ and the exceptional root system $G_{2}$. The embedding $G_{2}\to
B_{3}$ was already discussed. 
 The dimension of the adjoint representation for the series $B_{r}$
is $2r^{2}+r$ and of quasi-minuscule representation it is $2r+1$, so the dimension of the algebra
$\gf$ should be $(r+1)(2r+1)$ with the $\mathrm{rank}\gf=r+1$. The only solution is $\gf=D_{r+1}$
and we have already seen that this case corresponds to strongly multiplicity free module and
Gelfand-Tsetlin basis.

In case (3) the diagram of the representation contains weights with the length smaller than the
length of short root, so it cannot be a projection of $\Delta_{\gf}$ to an integer subalgebra. 

We need to consider representations of class (4) only for $m=1,2$ since the root system $\Delta_{\sfr}$
can have roots of at most two different lengths. In the case $m=1$ we get strongly multiplicity free
modules of $A_{r}$ which were discussed above. And for algebra $A_{2}$ and $m=2$ the weight diagram
of representation $L^{2\omega_{1}}_{A_{2}}$ coincides with the root system of exceptional Lie algebra
$G_{2}$. This case appears in the special embedding $A_{2}\to B_{3}$. But for $A_{r}, r>2$ there are
no root systems that coincide with weight diagram of $L^{2\omega_{1}}_{A_{r}}$. 

So far we've considered all the cases when $\mathrm{rank}\gf-\mathrm{rank}\af=1$. To complete our
classification of splints we need the classification of all the representations where multiplicities
of all non-zero weights are equal to one. Fortunately, such a classification was obtained in
\cite{plotkin1998visual}. It consists of multiplicity-free representations and adjoint
representations of algebras $A_{r}, B_{r}, C_{r}$, $F_{4}$ and $G_{2}$. 

In order to have the embedding $\af\to \gf$ such that the adjoint representation of $\gf$ is
decomposed into two adjoint representations of $\af$ (and possibly several trivial
representations) we need to satisfy following conditions:
\begin{itemize}
\item Denote by $n_{\gf}$ the total number of roots in $\Delta_{\gf}$ and the projection $\Delta'$.
  Then
  \begin{equation}
    \label{eq:5}
    n_{\gf}=2n_{\af}.
  \end{equation}
\item Denote rank of $\gf$ by shorthand notation $r_{\gf}$. Then the dimension
  \begin{equation}
    \label{eq:6}
    \mathrm{dim}\gf=n_{\gf}+r_{\gf}\geq 2\mathrm{dim}\af=2n_{\af}+2r_{\af}
  \end{equation}

\end{itemize}


For the adjoint representation of algebra $A_{r}$ we see that the projection $\Delta'$ of the root system
$\Delta_{\gf}$ should consist of roots of subalgebra $\af=A_{r}$ with multiplicity 2. The total
number of roots for $A_{r}$ is $r(r+1)$, so $\Delta'$ and $\Delta_{\gf}$ contain $2r(r+1)$ roots,
since there are no roots orthogonal to $\Delta_{\af}$. There are two possible solutions for
$\gf$: $D_{r+1}$ and $G_{2}$ for $r=2$. But the dimension of $D_{r+1}$ is $(r+1)(2r+1)$ while the
dimension of $A_{r}$ is $r(r+2)$ and $2\mathrm{dim}A_{r}>\mathrm{dim}D_{r+1}$, so the adjoint
representation of $D_{r+1}$ cannot be decomposed as twice the adjoint representation of $A_{r}$.
There is no corresponding embedding $A_{r}\to D_{r+1}$ and no splint. We see that conditions
\eqref{eq:5}\eqref{eq:6} do not hold. 

Let's consider the embedding $B_{r}\to \gf$ such that the adjoint representation of $\gf$ is
decomposed into two adjoint representations of $B_{r}$. It is straightforward to check that there is
no solutions for $\gf$ satisfying the \eqref{eq:5}\eqref{eq:6} for the series $A$, $B$, $C$,
$D$ and all exceptional simple Lie algebras. 

Since the number of roots of root system $C_{r}$ is the same as for $B_{r}$ and we need to check the
same cases, we see that there is no embedding for $C_{r}$ too. 

The adjoint representation of $G_{2}$ gives us the algebra $D_{4}$ as the solution of the
constraints \eqref{eq:5}\eqref{eq:6}. The embedding $G_{2}\to D_{4}$ exists, since there are
embeddings $G_{2}\to B_{3}$ and $B_{3}\to D_{4}$, but the decomposition of adjoint representation of
$D_{4}$ is different: $ \mathrm{ad}_{D_{4}}=\mathrm{ad}_{G_{2}}\oplus 2 L^{\omega_{1}}_{G_{2}}$. So
this case does not produce splint for the projection $\Delta'$. 

For the adjoint representation of $F_{4}$ there is no solution that satisfies conditions
\eqref{eq:5}\eqref{eq:6}.

The complete classification of splints for special embeddings $\Delta'=\pia(\Delta_{\gf})
\equiv \Delta_{\af}\cup \Delta_{\sfr}$ where $\gf, \af$ are simple and $\sfr$ is semisimple is:
\begin{itemize}
\item $B_{n}\to D_{n+1}$
\item $G_{2}\to B_{3}$
\end{itemize}

Note that in the splints for special embeddings there is no case where $\Delta_{\sfr}$ is embedded
non-metrically which is a direct result of our exhaustive classification of suitable characteristic
representations of the subalgebra $\af$.

Having obtained the complete classification of splints of the projections of
the root systems for special embeddings, we need to check whether such splints lead
to the coincidence of branching coefficients with weight multiplicities in modules of
the algebra $\sfr$.

Let's consider the first case - the embedding $B_{n}\to D_{n+1}$. The root system of $D_{n+1}$ has
exactly two simple roots that are different from the simple roots of $B_{n}$. In the notation of
\cite{bourbaki2002lie} they are $\alpha_{n}=e_n-e_{n+1}$ and $\alpha_{n+1}=e_n+e_{n+1}$ while
$\alpha_n^{B_n}=e_n$. This difference leads to the crucial difference in the Weyl groups of the
algebras: while $W_{B_{n}}$ is a semidirect product of the group of permutation $e_{i} \to e_{j}$
and the group of the change of sign $e_{i}\to (\pm 1)_{i}e_{i}$, $W_{D_{n+1}}$ is the same but for
the additional condition on the group of the change of sign: $\prod_{i} (\pm)_{i}=1.$ Since the
projection of $D_{n+1}$ on $B_{n}$ acts as
$(e_1,e_2,\dots,e_{n},e_{n+1})\to(e_1,e_2,\dots,e_{n},0)$, the singular element
$\Phi^{(\mu)}_{D_{n+1}}$ which is the orbit of $W_{D_{n+1}}$ will become a composition of the orbits
of $W_{B_{n}}$ after the projection.

Indeed, if $\mu+\rho_{D_{n+1}}=(a_1,a_2,\dots,a_{n+1})$ in the standart basis
$(e_1,e_2,\dots,e_{n+1})$, then the Weyl group $W_{D_{n+1}}$ will act on it by permutating $a_i$ and
changing their signs. After the projection the last coefficient in every element will be cut. Among
these elements will be groups in which all elements will have the same set of $a_i$ standing in
arbitrary order and having an arbitrary sign. These groups will be similar to the orbits of the Weyl
group $W_{B_{n}}$ but for the singular miltiplicity. A singular miltiplicity is a determinant
$\varepsilon (w)$ of the element $w$ of the Weyl group. Because of the additional condition on the
Weyl group $W_{D_{n+1}}$ the change of the signs of $a_i$ doesn't change the singular multiplicity.
That's not true for the Weyl group $W_{B_{n}}$, thus exactly half of the elements of each newly
formed orbit $\tilde\Phi^{\tilde\nu}_{B_n}$ of $W_{B_{n}}$ will have the wrong singular
multiplicity. Moreover, these orbits can be rewritten in the following form:

\begin{equation}
\tilde\Phi^{\tilde\nu}_{B_n}=\sum_{w\in W_{D_n}} \varepsilon(w) (1+s_{e_n})w e^{\tilde\nu+\rho_{B_n}}.
\end{equation}

This decomposition is possible because the orbit of the Weyl group of $D_n$ coincides with the half
of $\tilde\Phi^{\tilde\nu}_{B_n}$ that has correct singular multiplicities and other half can be
obtained by the action of the element $s_{e_n}$ of $W_{B_{n}}$.

\begin{figure}[h]
\begin{center}
\begin{minipage}[h]{0.5\linewidth}
\includegraphics[width=0.95\linewidth]{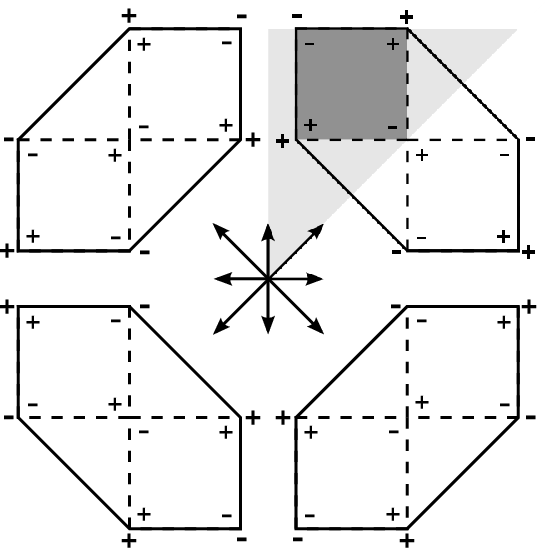}
\caption{Splint $D_3= B_2\cup A_1+A_1$: projection of the singular element $\Phi_{D_3}^{[1,1,2]}$ can be made from singular elements $\Phi_{A_1+A_1}^{[1,1]}$. Light-grey area is the fundamental Weyl chamber of subalgebra $B_2$.} 
\label{ris1} 
\end{minipage}
\hfill
\begin{minipage}[h]{0.47\linewidth}
\includegraphics[width=0.9\linewidth,viewport=0 0 157 157,clip]{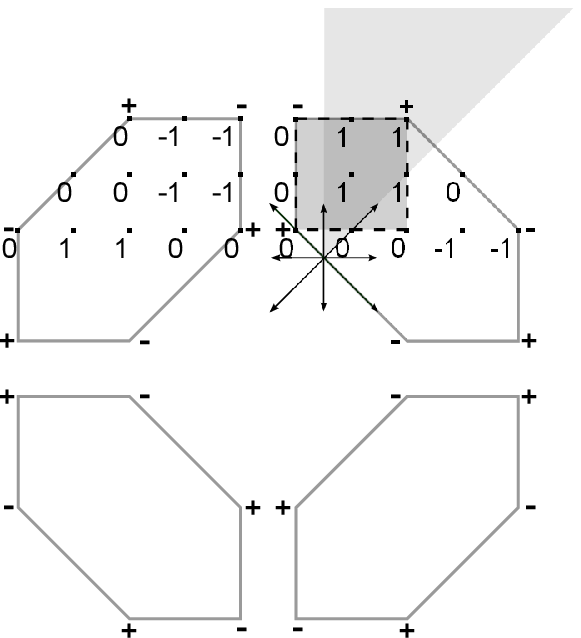}
\caption{Injection fan applied to the singular element $\Phi_{A_1+A_1}^{[1,1]}$ yields branching coefficients that are equal to 1.}
\label{ris2}
\end{minipage}
\end{center}
\end{figure}

Thus, the projection $(\Phi^{(\mu)}_{D_{n+1}})'$ consists of $(n+1)$ quasi-orbits (with ``wrong'' signs) of the Weyl group
$W_{B_{n}}$. It means that there are $(n+1)$ weights in the fundamental Weyl chamber $\bar C_{B_n}$
of $B_n$. It can be easily shown that these weights are $(a_1,a_2,\dots,a_{n})$ and weights obtained
from $(a_1,a_2,\dots,a_{n})$ by consequent subtraction of $\mu_i e_i$ starting with $\mu_n e_n$,
where $\mu_i$ is Dynkin labels of $\mu$ plus $1$. As it turns out, if one were to add to this set
the required weights with corresponding singular multiplicities to construct the singular element of
$\sfr=A_1+A_1+\dots+A_1$ as was done in the Introduction of this paper, all additional weight would
lie on the boundaries of the fundamental Weyl chamber $\bar C_{B_n}$. This fact is easily proved by
observing that the action of one of $s_{\alpha_i^{B_n}}$ on the weights doesn't change those
weigths. Thus, the projection $(\Phi^{\mu}_{D_{n+1}})'$ can be viewed as a quasi-orbit of the
singular element of $\sfr$:

\begin{equation}
\Phi^{[\mu_1,\mu_2,\dots,\mu_{n+1}]}_{D_{n+1}}=\sum_{w\in W_{D_n}} \varepsilon(w) (1+s_{e_n})w (e^{\mu'+\rho_{D_{n+1}}'-\tilde\mu}\Psi^{[\mu_1,\mu_2,\dots,\mu_{n}]}_{\sfr}),
\end{equation}
where $\tilde\mu=[\mu_1,\mu_2,\dots,\mu_{n}]$. This decomposition yields the branching rules similar to (\ref{bran1}):

\begin{equation}
b_{\left( \mu -\phi \left( \widetilde{\mu }-\widetilde{\nu }\right) \right)
}^{[\mu_1,\mu_2,\dots,\mu_{n+1}]}=M_{\left( \sfr\right) \widetilde{\nu }}^{[\mu_1,\mu_2,\dots,\mu_{n}]}, 
\label{bran1}
\end{equation}
which is in total agreement with the Gelfand-Tsetlin branching rules as the weights of the representations of $\sfr=A_1+A_1+\dots+A_1$ are always equal to one (see Fig.\ref{ris1} and Fig.\ref{ris2}).

Careful examination of the embedding $G_{2}\to B_{3}$ leads to exactly the  same conclusion. 

So far we only considered special splints where both embedded systems are root systems. However, it
is possible to study the splints where it is not the case. Such splints obviously can not serve to
simplify calculation of the branching rules but one may use the properties of special embeddings and
the projections to creat a perepesentation theory for the systems that are not root systems.

As the result we see that branching coefficients coincide with weight multiplicities for the splints
in the following table:
\begin{equation}
\label{eq:1}
\begin{array}{cc||c|c}
\hbox{type} & \hspace{0.25in}\Delta \hspace{0.25in} & \hspace{0.25in}\Delta
_{\frak{a}}\hspace{0.25in} & \hspace{0.25in}\Delta _{\sfr}\hspace{0.25in}
\\ \hline\hline
\hbox{(i)} & G_{2} & A_{2} & A_{2} \\
& F_{4} & D_{4} & D_{4} \\ 
\hbox{(*)} & B_{3} & G_{2} & A_{2}  \\
\hbox{(*)} & D_{r+1} & B_{r} & \oplus ^{r}A_{1}  \\
\hline
\hbox{(ii)} & B_{r}(r\geq 2) & D_{r} & \oplus ^{r}A_{1} \\
\hbox{(iii)} & A_{r}(r\geq 2) & A_{r-1}\oplus u\left( 1\right)  & \oplus
^{r}A_{1} \\
& B_{2} & A_{1}\oplus u\left( 1\right)  & A_{2}
\end{array}
\end{equation} 
here the splints marked with (*) are special splints.

\section*{Conclusion}
\label{sec:conclusion}

Computation of branching coefficients is important for different physical models with symmetry
breaking. This computation is drastically simplified if branching coefficients coincide with weight
multiplicities, since efficient Freudental formula can be used \cite{moody1982fast}. Classification
of splints for regular and special embeddings gives us all the cases when this coincidence takes
place. Aside from computational importance, this coincidence is very interesting from
representation-theoretic point of view, because it follows from new unexpected connection between
different simple Lie algebras. 

There exist special embeddings of semisimple Lie-algebras into simple Lie-algebras. For such
embeddings the procedure of finding splints that allow to calculate branching coefficient is similar
to the one conducted above. While the results of the procedure are not presented in this paper they
can be obtained with little difficulty.

\bibliography{special-bibliography}{} 
\bibliographystyle{apalike}
\end{document}